\documentclass[11pt]{amsart2000}
\input diagrams.tex
\newtheorem{pro}{Proposition}[section]
\newtheorem{teo}{Theorem}[section]
\newtheorem{cor}{Corollary}[section]
\newtheorem{lem}{Lemma}[section]
\newtheorem*{smt}{Split Monadicity Theorem under presence of coequalizers}
\theoremstyle{definition}
\newtheorem{exs}{Examples}[section]
\newtheorem{obs}{Remark}[section]
\newtheorem{defi}{Definition}[section]

\def\bkK{\mathbb K}
\def\bkT{\mathbb T}
\def\A{\mathcal{A}}
\def\B{\mathcal{B}}
\def\C{\mathcal{C}}
\def\E{\mathcal{E}}
\def\M{\mathcal{M}}
\def\N{\mathcal{N}}
\def\X{\mathcal{X}}
\def\S{\mathcal{S}et}
\def\T{\mathcal{T}\!op}
\def\Inj{\mbox{\scriptsize Inj}}
\def\reg{\mbox{\scriptsize reg}}

\def\op{\mbox{\scriptsize op}}

\begin{document}
\setlength{\unitlength}{0.01in}
\linethickness{0.01in}
\begin{center}
\begin{picture}(474,66)(0,0) 
\multiput(0,66)(1,0){40}{\line(0,-1){24}}
\multiput(43,65)(1,-1){24}{\line(0,-1){40}}
\multiput(1,39)(1,-1){40}{\line(1,0){24}}
\multiput(70,2)(1,1){24}{\line(0,1){40}}
\multiput(72,0)(1,1){24}{\line(1,0){40}}
\multiput(97,66)(1,0){40}{\line(0,-1){40}} 
\put(143,66){\makebox(0,0)[tl]{\footnotesize Proceedings of the Ninth Prague Topological Symposium}}
\put(143,50){\makebox(0,0)[tl]{\footnotesize Contributed papers from the symposium held in}}
\put(143,34){\makebox(0,0)[tl]{\footnotesize Prague, Czech Republic, August 19--25, 2001}}
\end{picture}
\end{center}
\vspace{0.25in}
\setcounter{page}{295}
\title[Pushout stability of embeddings]{Pushout stability of embeddings,
injectivity and categories of algebras}
\author{Lurdes Sousa}
\thanks{Partial financial support by Centro de Matem\'atica da 
Universidade de Coimbra and by Escola Superior de Tecnologia de Viseu is 
acknowledged}
\address{Departamento de Matem\'{a}tica\\
Escola Superior de Tecnologia\\
Instituto Polit\'{e}cnico de Viseu\\
Repeses 3500 Viseu, Portugal}
\email{sousa@mat.estv.ipv.pt}
\thanks{This article will be revised and submitted for publication 
elsewhere.}
\thanks{Lurdes Sousa,
{\em Pushout stability of embeddings, injectivity and categories of 
algebras},
Proceedings of the Ninth Prague Topological Symposium, (Prague, 2001),
pp.~295--308, Topology Atlas, Toronto, 2002}
\begin{abstract}
In several familiar subcategories of the category $\T$ of topological
spaces and continuous maps, embeddings are not pushout-stable. 
But, an interesting feature, capturable in many categories, namely in 
categories $\mathcal{B}$ of topological spaces, is the following: 
For $\mathcal{M}$ the class of all embeddings, the subclass of all
pushout-stable $\mathcal{M}$-morphisms (that is, of those
$\mathcal{M}$-morphisms whose pushout along an
arbitrary morphism always belongs to $\mathcal{M}$) is of the form
$A^{\mbox{\Inj}}$ for some space $A$, where $A^{\mbox{\Inj}}$
consists of all morphisms $m:X\rightarrow Y$ such that the map
Hom$(m,A): \mbox{Hom}(Y,A)\rightarrow \mbox{Hom}(X,A)$ is
surjective. 
We study this phenomenon. We show that, under mild assumptions, the
reflective hull of such a space $A$ is the smallest
$\mathcal{M}$-reflective subcategory of $\mathcal{B}$; furthermore, the
opposite category of this reflective hull is equivalent to a reflective
subcategory of the Eilenberg-Moore category $\S^\bkT$, where $\bkT$ is
the monad induced by the right adjoint Hom$(-,A):\T^{\op}\rightarrow \S$. 
We also find conditions on a category $\mathcal{B}$ under which the
pushout-stable $\mathcal{M}$-morphisms are of the form
$\mathcal{A}^{\mbox{\Inj}}$ for some category $\mathcal{A}$.
\end{abstract}
\subjclass[2000]{18A20, 18A40, 18B30, 18G05, 54B30, 54C10, 54C25}
\keywords{embeddings, injectivity, pushout-stability,(epi)reflective
subcategories of $\T$, closure operator, Eilenberg-Moore categories}
\maketitle

\setcounter{section}{-1}

\section{Introduction}

It is well-known that in the category $\T$ of topological spaces
and continuous maps, the class of all embeddings is pushout-stable, that
is, if we have a pushout in $\T$,
$$
\begin{diagram}
\bullet&\rTo^{m}&\bullet\\
\dTo^{f} &&\dTo_{\overline{f}}\\
\bullet&\rTo^{\overline{m}}&\bullet\\
\end{diagram}
$$
with $m$ an embedding, then $\overline{m}$ is also an embedding.
But this is not true for several subcategories of $\T$, namely, it is not
true for various epireflective subcategories of $\T$.

Let $\B$ be a subcategory\footnote{Along this paper we assume that all
subcategories are full.} of $\T$, let $\mathcal{M}$ denote the class of
embeddings in $\B$.
By $P_\B(\M)$ we denote the class of all morphisms $m\in \M$ such that
every pushout in $\B$ of $m$ along any morphism belongs to $\M$.
In this paper, this class of morphisms is shown to have an important
r\^ole concerning some nice behaviours of the reflective hull of a
topological space.

In several reflective subcategories of $\T$, it is possible to select a
subcategory $\A$ such that $P_\mathcal{B}(\M)=\A^{{\Inj}_\B}$, where
$\A^{{\Inj}_\B}$ denotes the class of all morphisms $m$ of $\B$ such that
each $A\in \A$ is $m$-injective, that is, the map 
Hom$(m,A): \mbox{Hom}(Y,A)\rightarrow \mbox{Hom}(X,A)$ 
is surjective. 
Under smooth conditions, this fact implies that the reflective hull of
$\A$ is just the smallest $\M$-reflective subcategory of the epireflective
hull of $\A$. 
We recall that an $\M$-reflective subcategory is a reflective subcategory
whose corresponding reflections belong to $\M$. 
The most interesting case, which is rather common, is when $\A$ is chosen
as one-object subcategory. 
In fact, if $A$ is a topological space for which $P_\mathcal{B}(\M)$
consists of all morphisms $m$ such that $A$ is $m$-injective, then, it
often holds that the dual category of the reflective hull of $A$ in $\T$
is a reflective subcategory of $\S^\bkT$ for $\bkT$ the monad induced by
Hom$(-,A)$. 
Furthermore, under convenient requisites, it is equivalent to $\S^\bkT$.

As very helpful tools, two Dikranjan-Giuli closure operators will be used: 
the regular closure operator (introduced in \cite{Sa,DG}) and the
orthogonal closure operator (introduced in \cite{S1}). 
We shall see that the equality $P_\mathcal{B}(\M)=\A^{{\Inj}_\B}$, 
combined with mild conditions, implies that the two closure operators
coincide in $P_\mathcal{B}(\M)$, and this plays an important r\^ole in the
main result of this paper, Theorem \ref{theo}.

\section{Injectivity and pushout-stability}

Let $A$ be a topological space and let $f:X\rightarrow Y$ be a morphism
in $\T$.
We say that $A$ is {\it $f$-injective} (respectively, {\it orthogonal to 
$f$})
whenever the map
$\mbox{Hom}(f,A):\mbox{Hom}(Y,A)\rightarrow \mbox{Hom}(X,A)$ is
surjective (respectively, bijective).
For $\A$ a class of spaces, that is, a full subcategory of $\T$, we denote
by $\A^{\Inj }$ the class of those morphisms $f$ such that all objects of
$\A$ are $f$-injective. 
In case $\A$ has just one object $A$, we write $A^{\Inj}$. 
Analogously, $\A^\bot$ is the class of morphisms $f$ such that any object
of $\A$ is orthogonal to $f$. 
Given a class $\N$ of morphisms, the subcategory of $n$-injective objects
(respectively, objects orthogonal to $n$), for all $n\in \N$, is named
$\N_{\Inj}$ (respectively, $\N_\bot$).

Let $\B$ be a subcategory of $\T$. 
We denote by $\A^{{\Inj}_\B}$ the class $\A^{\Inj}$ restricted to
$\B$-morphisms. 
$P_\B(\M)$ designates the class of all $\M$-morphisms whose pushout in
$\B$ along any morphism exists and belongs to $\M$.

Throughout this paper, unless something is said on the contrary, $\M$
denotes the class of embeddings of $\T$.

An interesting feature, capturable in various reflective subcategories
$\B$ of $\T$, is that the class $P_{\B}(\M)$ coincides with
$\A^{{\Inj}_\B}$ for some subcategory $\A$, in many cases this subcategory
$\A$ having just one object. 
Next, some examples of this occurrence are given.

\begin{exs}\label{Ex1.1}
\mbox{}
\begin{enumerate}
\item
Since in $\T$ embeddings are pushout-stable, we have the equality 
$P_{\T}(\M)=\M$. 
Moreover, $\M=A^{\Inj}$ for $A$ the topological space $\{0,1,2\}$ whose
only non trivial open is $\{ 0\}$.
\item
For the subcategory ${\T}_{\, 0}$ of $T_0$-spaces, again 
$P_{{\T}_{\, 0}}(\M)=\M$. 
In this case, it holds $\M=S^{{\Inj}_{\T_0}}$, where $S$ is the 
Sierpi\'nski space.
\item
For $\B$, the subcategory of 0-dimensional spaces, $\M \not= P_{\B}(\M)$,
but again $P_{\B}(\M)=A^{{\Inj}_{\B}}$, where $A$ is the space $\{0,1,2\}$
whose topology has as only non trivial opens $\{0\}$ and $\{ 1,2\}$. 
The morphisms of $P_\B(\M)$ are just those embeddings $m:X\rightarrow Y$
such that for each clopen set $G$ of $X$ there is some clopen $H$ in $Y$
such that $G=m^{-1}(H)$ (cf.\ \cite{S2}).
\item
If $\B$ is the subcategory of 0-dimensional Hausdorff spaces, 
$\M \not= P_{\B}(\M)=D^{{\Inj}_{\B}}$, where $D$ is the space $\{0,1\}$
with the discrete topology. 
Analogously to the above example, the morphisms of $P_\B(\M)$ are those
embeddings $m:X\rightarrow Y$ in $\B$ such that each clopen in $X$ is the
inverse image by $m$ of some clopen in $Y$ (cf.\ \cite{S2}).
\item
For $\mathcal{T}ych$ the subcategory of Tychonoff spaces, 
$\M \not= P_{\mathcal{T}ych}(\M)=I^{{\Inj}_{\mathcal{T}ych}}$, where $I$
is the unit interval. 
The $P_{\mathcal{T}ych}(\M)$-morphisms are just the $C^*$-embeddings
(cf.\ \cite{S2}).
\item
In the category $\mathcal{I}nd$ of indiscrete spaces embeddings are 
pushout-stable, and we have that 
$P_{\mathcal{I}nd}(\M) = \M = \{0,1\}^{{\Inj}_{\mathcal{I}nd}}$.
\end{enumerate}
\end{exs}

\begin{exs}
This phenomenon arises in several other categories, not just reflective
subcategories of $\T$. 
Some examples are given in the following.
\begin{enumerate}
\item
For $\mathcal{A}b$ the category of abelian groups and homomorphisms of
group, and $\M$ the class of all monomorphisms, $P(\M)=\M=\A^{\Inj}$,
where $\A$ is the subcategory of divisible abelian groups. 
Similarly, in the category of torsion-free abelian groups and
homomorphisms of group, and $\M$ the class of all monomorphisms,
$P(\M)=\M=\A^{\Inj}$, where $\A$ is the subcategory of divisible groups.
\item
For $\mathcal{V}ec_{\bkK}$ the category of linear spaces and linear maps
over $\bkK$, and $\M$ the class of all monomorphisms, 
$P(\M)=\M=\bkK^{{\Inj}}$.
\item 
Recall that a {\it separated quasi-metric space} is a pair $(X,d)$, where
$X$ is a set and $d$ is a function with domain $X\times X$ which satisfies
the same axioms as a metric in $X$, but $d$ may take the value $\infty$. 
A separated quasi-metric space is said to be \textit{complete} if all its
Cauchy sequences converge. 
A function $f$ between two separated quasi-metric spaces $(X,d)$ and
$(Y,e)$ is \textit{non-expansive} provided that $e(f(x),f(y))\leq d(x,y)$
for all $(x,y)\in X\times X$. 
Let $\mathcal{SM}et$ denote the category of separated
quasi-metric spaces and non-expansive maps and let 
$\mathcal{CSM}et$ denote its subcategory of complete
separated quasi-metric spaces. 
Then, for $\M$ the class of all non-expansive monomorphisms, $\M$ is
pushout-stable, so it coincides with $P(\M)$, and it holds 
$\M=\mathcal{CSM}et^{\Inj}$.
\item
For the category ${\M}et$ of metric spaces and non-expansive maps
and $\M$ defined similarly as in 3, we have that
$$
P(\M) = \{m\in \M | 
\mbox{ the pushout of $m$ along any morphism exists}\}
$$
and it coincides with $\mathcal{CM}et^{\Inj}$, where $\mathcal{CM}et$
denotes the subcategory of complete metric spaces. 
An analogous situation holds for the category of normed spaces and
non-expansive maps and the subcategory of Banach spaces.
\end{enumerate}
\end{exs}

One first question which comes from these examples is the following: 
When is $P_\B(\M)$ of the form $\A^{\Inj_\B}$ for some subcategory $\A$ of
a given category, in particular, for $\A$ a subcategory of $\T$? 
The next proposition is a partial answer for this.

Let us recall that, if $\N$ is a class of morphisms in a category $\B$, 
containing all isomorphisms and being closed under composition with
isomorphisms, we say that $\B$ \textit{has enough $\N$-injectives}
provided that, for each $B\in \B$, there is an $\N$-morphism
$n:B\rightarrow A$ with the codomain $A$ in $\N_{\Inj}$. 
A morphism $n:B\rightarrow A$ of $\mathcal{N}$ is said to be {\it
$\N$-essential} whenever any composition $k\cdot n$ belongs to $\N$ only
if $k\in \N$. The category $\B$ is said to have {\it $\N$-injective hulls}
if, for each $B\in \B$, there is some $\N$-essential morphism
$n:B\rightarrow A$ with $A$ $\N$-injective.

The following easy lemma to prove will be useful in this paper.

\begin{lem}
For any epireflective subcategory $\B$ of $\T$, the class $P_\B(\M)$ has
the following properties:
\begin{itemize}
\item[(i)]
It is closed under composition.
\item[(ii)]
If $m\cdot n \in P_\B(\M)$ then $n\in P_\B(\M)$.
\end{itemize}
\end{lem}

\begin{pro}\label{pro1}
Let $\B$ be an epireflective subcategory of $\T$.
\begin{enumerate}
\item
If $\B$ has enough $P_\B(\M)$-injectives, then $P_\B(\M)=\A^{{\Inj}_\B}$
for some subcategory $\A$ of $\B$.
\item
If $P_\B(\M)=\A^{{\Inj}_\B}$ for some $\M$-reflective subcategory $\A$ of
$\B$, then $\B$ has $P_\B(\M)$-injective hulls, and 
$\left[P_\B(\M)\right]_{{\Inj}_\B}=\A$.
\end{enumerate}
\end{pro}

\begin{proof} 
The fact that $\B$ is an epireflective subcategory of $\T$ ensures that
$\B$ has pushouts and an (Epi, Embedding)-factorization system for
morphisms.

(1). 
Let $\A=\left[P_\B(\M)\right]_{{\Inj}_\B}$; it is clear that 
$P_\B(\M)\subseteq \A^{{\Inj}_\B}$. 
In order to show the converse inclusion, let $f:X\rightarrow Y$ belong to
$\A^{{\Inj}_\B}$ and let $f=me$ be the (Epi, Embedding)-factorization of
$f$. 
By hypothesis, there is a $P_\B(\M)$-morphism $n:X\rightarrow A$ with
$A\in \left[ P_\B(\M)\right]_{{\Inj}_\B}$. 
Let $\overline{n}:Y\rightarrow A$ be such that $\overline{n}f=n$; then the
equality $n1_X=(\overline{n}m)e$ guarantees the existensce of a unique 
morphism $t$ such that $te=1_X$ and, consequently, $e$ is an isomorphism
and $f\in \M$. 
Now, to conclude that $f:X\rightarrow Y$ belongs to $P_\B(\M)$, let
$$
\begin{diagram}
X&\rTo^{f}&Y\\
\dTo^{g} &&\dTo_{\overline{g}}\\
Z&\rTo^{\overline{f}}&W\\
\end{diagram}
$$
be a pushout in $\B$. 
Since $\B$ has enough $P_\B(\M)$-injectives, there is some 
$P_\B(\M)$-morphism $p:Z\rightarrow A$, with $A\in \B$ being
$P_\B(\M)$-injective. 
Let $h:Y\rightarrow A$ be such that $hf=pg$; then, there is $t$ such that
$t\overline{f}=p$ and so the fact that $p\in \M$ implies that also
$\overline{f}\in \M$.

(2). 
Let us first show that each reflection $r_B: B\rightarrow RB$ of $B\in \B$
in $\A$ belongs to $P_\B(\M)$. 
Let $(\overline{f}, \overline{r})$ be the pushout of $(r_B, f)$ for some
morphism $f:B\rightarrow C$.
Since $Rf\cdot r_B=r_C\cdot f$, there is a unique $t$ such that 
$t\cdot \overline{r}=r_c$ and $t\cdot \overline{f}=Rf$. 
By hypothesis, $r_C\in \M$, and this implies that $\overline{r}\in \M$. 
Then $r_B\in P_\B(\M)$.

To show that $r_B$ is $P_\B(\M)$-essential, let $k:RB\rightarrow K$ be a
morphism such that $k\cdot r_B$ belongs to $P_\B(\M)$. 
It is easy to see that the fact that $\A$ is $\M$-reflective in $\B$
implies that each reflection in $\A$ is an epimorphism in $\B$.
Consequently, for any morphism $f:RB\rightarrow D$, $(1_D, f)$ is a
pushout of $(r_B, f\cdot r_B)$. 
Now, let $f:RB \rightarrow D$ be a morphism and let $(\overline{k},
\overline{f})$ be the pushout of $(k,f)$. 
Then, since $(1_D, f)$ is a pushout of $(r_B, f\cdot r_B)$, $\overline{k}$
is the pushout of $k\cdot r_B$ along $f\cdot r_B$.
Thus, $\overline{k}$ belongs to $\M$ and, therefore, $k\in P_\B(\M)$.

In order to conclude that $\left[P_\B(\M)\right]_{{\Inj}_\B}=\A$, since
the inclusion $\A\subseteq \left[P_\B(\M)\right]_{{\Inj}_\B}$ is clear, it
remains to show the reverse inclusion. 
Let $B\in \left[P_\B(\M)\right]_{{\Inj}_\B}$; the fact that the reflection
$r_B:B\rightarrow RB$ belongs to $P_\B(\M)$ implies that there is some
$t:RB\rightarrow B$ such that $t\cdot r_B=1_B$, and this implies that
$r_B$ is an isomorphism.
\end{proof}

Given a subcategory $\A$ of $\T$, let $\mathcal{R}(\A)$ denote the 
reflective hull of $\A$ in $\T$, provided it exists. 
It is well known that, when $\A$ is small, $\mathcal{R}(\A)$ exists and
coincides with the limit closure of $\A$.

The subcategories $\A$ such that $P_\B(\M)=\A^{{\Inj}_\B}$, for $\B$ its
epireflective hull, have a very special property, as established next.

\begin{pro}
Let $\B$ be the epireflective hull of $\A$ in $\T$. 
If $P_\B(\M)=\A^{{\Inj}_\B}$, then $\mathcal{R}(\A)$ is the smallest 
$\M$-reflective subcategory of $\B$ (provided $\mathcal{R}(\A)$ exists).
\end{pro}

\begin{proof}
Since $\B$ is the epireflective hull of $\A$, it is clear that
$\mathcal{R}(\A)$ is $\M$-reflective in $\B$ (cf.\ \cite{AHS}). 
Let $\C$ be another $\M$-reflective subcategory of $\B$. 
Then, as shown in the proof of 2 of Proposition \ref{pro1}, each
reflection $s_B:B\rightarrow SB$ of $B\in \B$ into $\C$ belongs to
$P_\B(\M)$. 
Therefore, given $B\in \A$, since 
$\A\subseteq \left[P_\B(\M)\right]_{{\Inj}_\B}$, there is some 
$t: SB\rightarrow B$ such that $t\cdot s_B=1_B$, and so $s_B$ is an 
isomorphism. 
Consequently, $\A\subseteq \C$ and, thus, $\mathcal{R}(\A)\subseteq \C$.
\end{proof}

\section{On closure operators}

In this section we recall some notions concerning closure operators, in 
the sense of Dikranjan-Giuli (\cite{DG,DGT,DT}), and obtain some results,
which are going to be useful in the next section.

Let $\N$ be a class of monomorphisms in a category $\X$, which contains
all isomorphisms, is closed under composition and is left-cancellable,
(i.e., if $m\cdot n$, $m\in \N$ then $n\in \N$). 
Considering $\N$ as a subcategory of the category $\X^2$ (of all morphisms
of $\X$), let $u:\N\rightarrow \X$ be the codomain functor, that is, the
functor which, to each morphism 
$(r,s):(m:X\rightarrow Y) \longrightarrow (n:Z\rightarrow W)$ 
of $\N$, assigns the morphism $s:Y\rightarrow W$ of $\X$.

A {\it closure operator} in $\X$ with respect to $\N$ consists of a
functor $c:\N\rightarrow \N$ such that $u\cdot c=u$ and of a natural
transformation $\delta :Id_\N\rightarrow c$ such that 
$u\cdot \delta=Id_u$. 
Therefore, a closure operator determines, for each $m:X\rightarrow Y$ in
$\N$, morphisms $c(m)$ and $d(m)$ and a commutative diagram
$$
\begin{diagram}
X&\rTo^{d(m)}&\overline{X}\\
\dTo^{m} &&\dTo_{c(m)}\\
Y&\rTo^{1_Y}&Y\\
\end{diagram}
$$
where $\delta_m=(d(m),1_Y):m\rightarrow c(m)$.

We recall that, if $\N$ is part of an $(\mathcal{E},\N)$-factorization
system for morphisms, then a closure operator $c:\N\rightarrow \N$ may be 
equivalently described by a family of functions
$$\left( c_X:\N_X\rightarrow \N_X \right)_{X\in \X},$$
where $c_X(m)=c(m)$ for each $m$, and $\N_X$ denotes the class of all
$\N$-morphisms with codomain in $X$, fulfilling the conditions:
\begin{itemize}
\item[(1)] 
$m\leq c_X(m)$,\ $m\in \N_X$;
\item[(2)]
if $m\leq n$, then $c_X(m)\leq c_X(n)$,\
$m,n \in \N_X$; and
\item[(3)]
$c_X(f^{-1}(m))\leq f^{-1}(c_Y(m))$, for each morphism $f:X\rightarrow Y$
and $m\in \N_Y$.
\end{itemize}

If $c$ is a closure operator in $\X$ with respect to $\N$, then a
morphism $(n:X\rightarrow Y)\in \N$ is said to be {\it $c$-dense} if
$c(n)\cong 1_Y$. 
It is said to be {\it $c$-closed} when $c(n)\cong n$. 
We say that $c$ is {\it weakly hereditary}, if $d(n)$ is $c$-dense for any
$n\in \N$, and {\it idempotent} if $c(n)$ is $c$-closed for any $n\in \N$. 
An object $X\in \X$ is said to be {\it absolutely $c$-closed} whenever
each $\N$-morphism with domain in $X$ is $c$-closed. 

\subsection*{The regular closure operator} 

Let $\A$ be a reflective subcategory of $\T$. 
The {\it regular closure operator} induced by $\A$, which will be denoted
by 
$c^{\mbox{\scriptsize reg}}_\A:\M\rightarrow\M$ (just
$c^{\mbox{\scriptsize reg}}$, without the index $\A$, when there is no
confusion about the subcategory $\A$ which is involved), assigns, to each
$m:X\rightarrow Y$ in $\M$, the equalizer eq$(u\cdot r_W,v\cdot r_W )$
where $u,v:Y\rightarrow W$ form the cokernel pair of $m$ and $r_W$ is a
reflection in $\A$.
The regular closure operator, which is an idempotent closure operator, in
the sense described above, has been widely studied, not only in the
topological context, but in a more general setting (see \cite{DG}).

It is easy to conclude that if $\B$ is the epireflective hull of $\A$ in
$\T$ (or, equivalently, $\A$ is $\M$-reflective in $\B$), then
$c^{\mbox{\scriptsize reg}}_\A=c^{\mbox{\scriptsize reg}}_\B$.

From now on, given a subcategory $\A$ of $\T$, $\E(\A)$ always denotes the
epireflective hull of $\A$ in $\T$.

\subsection*{The orthogonal closure operator} 

Let $\A$ be a reflective subcategory in its epireflective hull $\E(\A)$ in
$\T$. 
Let $r_X:X\rightarrow RX$ be the reflection of $X$ in $\A$. 
To each $m:X\rightarrow Y$ in $P_{\mathcal{E}(\A)}(\M)$, we assign the
$\M$-morphism $c^{\mbox{\scriptsize ort}}_Y(m):\overline{X}\rightarrow Y$
which is the pullback of the pushout, in $\E(\A)$, of $m$ along $r_X$,
according to the next diagram.
$$
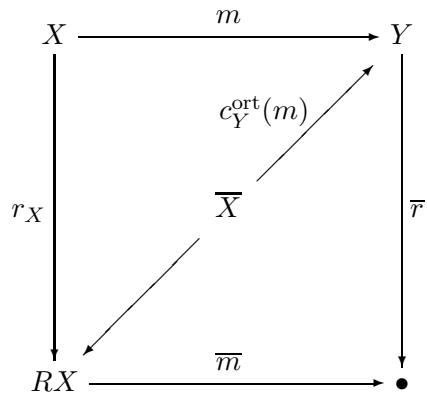
\begin{diagram}
X		&	&\rTo^{m}		&					&Y			\\
		&	&			&\ruTo{c^{\mbox{\scriptsize ort}}_Y(m)}	&			\\
\dTo^{r_X}	&	&\overline{X}		&					&\dTo_{\overline{r}}	\\
		&\ldTo{}&			&					&			\\
RX		&	&\rTo^{\overline{m}}	&					&\bullet		\\
\end{diagram}
$$

Using results of \cite{S1}, we get the next theorem. 

\begin{teo}\label{abs clo}
For every reflective subcategory $\A$ of $\T$, if 
$c^{\mbox{\scriptsize ort}}_Y(m):\overline{X}\rightarrow Y$ belongs to
$P_{\mathcal{E}(\A)}(\M)$ for all $m\in P_{\mathcal{E}(A)}(\M)$, then
$$c_\A^{\mbox{\scriptsize ort}}:P_{\mathcal{E}(\A)}(\M) \longrightarrow
P_{\mathcal{E}(\A)}(\M)$$
defines a weakly hereditary idempotent closure operator in
$\mathcal{E}(\A)$, and $\A$ is just the subcategory of all absolutely
$c^{\mbox{\scriptsize ort}}_\A$-closed objects.
\end{teo}

As expressed for $c^{\mbox{\scriptsize reg}}_\A$ before, we remove $\A$
from $c^{\mbox{\scriptsize ort}}_\A$ if there is no reason for ambiguity.

We are going to make use of the following definition.

\begin{defi}
We say that epimorphisms in a category $\B$ are \textit{left-cancell\-able
with respect to $\M$} provided that $q\in \mbox{Epi}(\B)$ whenever 
$p\cdot q \in \mbox{Epi}(\B)$ with $p,q\in \M$.
\end{defi}

\begin{cor}
If $\A$ is a reflective subcategory of $\T$ such that, in
$\mathcal{E}(\A)$, epimorphisms are left-cancellable with respect to $\M$,
then
$$c_\A^{\mbox{\scriptsize ort}}:P_{\mathcal{E}(\A)}(\M) \longrightarrow
P_{\mathcal{E}(\A)}(\M)$$
defines a weakly hereditary idempotent closure operator and $\A$ is just
the subcategory of all absolutely 
$c^{\mbox{\scriptsize ort}}_{\A}$-closed objects.
\end{cor}

\begin{proof}
It follows from the above proposition and the fact that the
left-cancellability of epimorphisms in $\mathcal{E}(\A)$ with respect to
$\M$ implies that 
$c^{\mbox{\scriptsize ort}}_Y(m):\overline{X}\rightarrow Y$ 
belongs to $P_{\mathcal{E}(\A)}(\M)$ for all 
$m\in P_{\mathcal{E}(A)}(\M)$, as it is shown in Remark 4.4 of \cite{S1}.
\end{proof}

\begin{lem}\label{M}
In any reflective subcategory of $\T$ where embeddings are
pushout-stable, epimorphisms are left-cancellable with respect to
embeddings.
\end{lem}

\begin{proof} 
First, we point out the following fact: 
In any category with pushouts, a morphism $f$ is an epimorphism iff 
$(f,f)$ has a pushout of the form $(f^\prime, f^\prime)$ for some morphism
$f^\prime$; in this case, $f^\prime$ is an isomorphism.

Let $\B$ be a reflective subcategory of $\T$ and let $m:X\rightarrow Y$
and $n:Y\rightarrow Z$ be $\M$-morphisms of $\B$ such that $nm$ is an
epimorphism. 
Consider the following diagram, where the smallest four squares are
pushouts in $\B$.
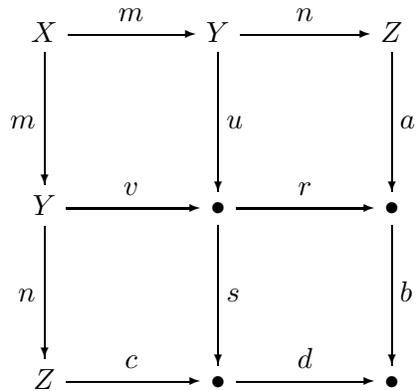
$$
\begin{diagram}
X	&\rTo^{m}	&Y		&\rTo^{n}	&Z		\\
\dTo^{m}&		&\dTo_{u}	&		&\dTo_{a}	\\
Y	&\rTo^{v}	&\bullet	&\rTo^{r}	&\bullet	\\
\dTo^{n}&		&\dTo_{s}	&		&\dTo_{b}	\\
Z	&\rTo^{c}	&\bullet&	\rTo^{d}	&\bullet	\\
\end{diagram}
$$
Then $b\cdot a=d\cdot c$ and we have that
$$(br)u=ban=dcn=dsv=(br)v.$$
Since $br$ is the composition of two embeddings, it is an embedding, and
thus $u=v$; consequently, $m$ is an epimorphism.
\end{proof}

\begin{obs}
The above lemma provides several examples of subcategories of $\T$
where epimorphisms are left-cancellable with respect to embeddings. 
This occurrence is far from to confine to subcategories under the
condition of Lemma \ref{M}. 
In fact, it is a very common property in epireflective subcategories of
$\T$.
This is the case, for instance, for $T_0$-spaces, $T_1$-spaces, Hausdorff
spaces, completely regular Hausdorff spaces, compact Hausdorff spaces,
compact 0-dimensional Hausdorff spaces. 
The left-cancellability of epimorphisms with respect to embeddings, in 
these categories, is easy to conclude, taking into account that in each
one of them epimorphisms are the surjective continuous maps or the dense
continuous maps (cf.\ \cite{DT,KMPT}).
\end{obs}

We are going to make use of the following lemma.

\begin{lem}[cf.\ \cite{S1}]\label{2.1}
Let $\A$ be a subcategory of $\T$. 
Then, in $\mathcal{E}(\A)$,
$$\begin{array}{lll}
\A^\bot&
=&
\left( \left(\A^\bot \right)_\bot \right)^\bot\\ 
&
=&
\{
c_\A^{\mbox{\scriptsize ort}}\mbox{-dense } 
P_{\mathcal{E}(\A)}\mbox{-morphisms}
\}\\
&
\subseteq&
P_{\mathcal{E}(\A)}(\M) \cap \mbox{Epi}(\mathcal{E}(\A)).
\end{array}$$
\end{lem}

We point out that, although the existence of the reflective hull of any
subcategory of $\T$ cannot be guaranteed (\cite{TAR}), we can state the
following:

\begin{pro}\label{Kelly}
Let $\A$ be a subcategory of $\T$. If $\mathcal{R}(\A)$ exists, then it
coincides with $\left(\A^\bot \right)_\bot$.
\end{pro}

\begin{proof} 
From Proposition 3.1.2 and Theorem 4.1.3 of \cite{FK}, one derives that if
$\X$ is a category fulfilling the following conditions:
\begin{itemize}
\item
it is complete, cocomplete and co-wellpowered;
\item
it has an $(\mathcal{E}, \M)$-factorization system for morphisms
with $\mathcal{E} = \mbox{Epi}(\X)$;
\item
it has a generator;
\item
for every countable family $(m_i:C_i\rightarrow B)_{i\in \omega}$
of $\M$-subobjects of $B\in\X$, and any morphism $g$ with codomain in
$B$, it holds the equality 
$\bigvee_{i\in \omega}g^{-1}(m_i)=g^{-1}(\bigvee_{i\in \omega} m_i)$;
\end{itemize}
then, for each morphism $f$ of $\X$, the subcategory $\{f\}_\bot$ is
reflective. 
The category $\T$ is under these conditions for $\M$ the class of all
embeddings. 
Consequently, given a subcategory $\A$ of $\T$ which has a reflective
hull, we have the inclusions:
$$
\begin{array}{lll}
\mathcal{R}(\A)&
\subseteq&
\bigcap_{f\in \A^\bot}\{f\}_\bot\\
&
=&
\left(\bigcup_{f\in \A^\bot}\{ f\}\right)_\bot\\
&
=&
\left(\A^\bot \right)_\bot.
\end{array}
$$
Therefore, since the inclusion 
$\left(\A^\bot \right)_\bot\subseteq \mathcal{R}(\A)$
always occurs, we conclude that 
$\mathcal{R}(\A)=\left(\A^\bot \right)_\bot $.
\end{proof}

The next proposition establishes conditions under which the closure
operators 
$c^{\mbox{\scriptsize ort}}_\A$ and 
$c^{\mbox{\scriptsize reg}}_\A$ coincide in $P_{\mathcal{E}(\A)}(\M)$. 
This result is going to be very useful in the proof of the theorem of the
next section.

\begin{pro}\label{r=c}
Let $\A$ be a subcategory of $\T$ with reflective hull.
If epimorphisms in $\mathcal{E}(\A)$ are left-cancellable with respect to
$\M$, 
$c_\A^{\mbox{\scriptsize reg}}$ is weakly hereditary in $\mathcal{E}(\A)$ and 
$P_{\mathcal{E}(\A)}(\M)=\A^{{\Inj\,}_{\mathcal{E}(\A)}}$, then
$$
c_{\mathcal{R}( \A)}^{\mbox{\scriptsize reg}} = 
c_{\mathcal{R}(\A)}^{\mbox{\scriptsize ort}}\
\text{in}\ P_{\mathcal{E}(\A)}(\M).$$
\end{pro}

\begin{proof} 
In order to simplify the writing, let $\A^{\mbox{\scriptsize Inj}}$ and
$\A^\bot$ refer to just morphisms in $\mathcal{E}(\A)$. 
It is easy to conclude that the fact that $\mathcal{E}(\A)$ is the
epireflective hull of $\A$ implies that 
$$\A^\bot \subseteq 
P_{\mathcal{E}(\A)}(\M)\bigcap \mbox{Epi}(\mathcal{E}(\A))$$ 
(see Lemma \ref{2.1}). 
On the other hand, it is clear that 
$$\A^{\mbox{\scriptsize Inj}}\bigcap \mbox{Epi}(\mathcal{E}(\A)) \subseteq
\A^\bot.$$
Therefore, we have that
$$
\begin{array}{lll}
\A^\bot&
\subseteq&
P_{\mathcal{E}(\A)}(\M)\bigcap \mbox{Epi}(\mathcal{E}(\A))\\
&
=&
\A^{\mbox{\scriptsize Inj}}\bigcap \mbox{Epi}(\mathcal{E}( \A))\\
&
\subseteq&
\A^\bot,
\end{array}
$$
and so,
\begin{equation}\label{ortogonal}
\A^\bot = 
P_{\mathcal{E}(\A)}(\M)\bigcap \mbox{Epi}(\mathcal{E}(\A)).
\end{equation} 
Furthermore, from Lemma \ref{2.1} and Proposition \ref{Kelly},
\begin{equation}\label{A^t=(R(A))^t}
\begin{array}{lll}
\A^\bot&
=&
\left( (\A^\bot)_\bot\right)^\bot\\
&
=&
\left(\mathcal{R}(\A)\right)^\bot\\ 
&
=& 
\{c_{\mathcal{R}(\A)}^{\mbox{\scriptsize ort}}\mbox{-dense } 
P_{\mathcal{E}(\A)}\mbox{-morphisms}\}.
\end{array}
\end{equation} 
It is well-known that the 
$c_{\mathcal{R}(\A)}^{\mbox{\scriptsize reg}}$-dense morphisms are just
the epimorphisms of $\mathcal{E}(\A)$ (see, for instance, \cite{DT}).
Consequently, using (\ref{ortogonal}) and (\ref{A^t=(R(A))^t}), we get
that a morphism of $P_{\mathcal{E}(\A)}(\M)$ is 
$c_{\mathcal{R}(\A)}^{\mbox{\scriptsize reg}}$-dense iff it is 
$c_{\mathcal{R}(\A)}^{\mbox{\scriptsize ort}}$-dense.

Let now $m$ be a morphism in $P_{\mathcal{E}(\A)}(\M)$.
Since the orthogonal closure operator is always smaller than the regular
one (cf.\ \cite{S1}), there is some morphism $d$ such that 
$c_{\mathcal{R}(\A)}^{\mbox{\scriptsize ort}}(m) = 
c_{\mathcal{R}(\A)}^{\mbox{\scriptsize reg}}(m)\cdot d$. 
Let $d^\prime$ be the morphism such that
$m = c_{\mathcal{R}(\A)}^{\mbox{\scriptsize ort}}(m)\cdot d^\prime$. 
Since 
$c_{\mathcal{R}(\A)}^{\mbox{\scriptsize reg}}$ is weakly hereditary, the 
morphism $d\cdot d^\prime$ is 
$c_{\mathcal{R}(\A)}^{\mbox{\scriptsize reg}}$-dense, and, thus, also
$c_{\mathcal{R}(\A)}^{\mbox{\scriptsize ort}}$-dense. 
Since 
$c_{\mathcal{R}(\A)}^{\mbox{\scriptsize ort}}$ is a weakly hereditary
idempotent closure operator, it follows that it has a
($c_{\mathcal{R}(\A)}^{\mbox{\scriptsize reg}}$-dense,
$c_{\mathcal{R}(\A)}^{\mbox{\scriptsize reg}}$-closed)-factorization
system with respect to $P_{\mathcal{E}(\A)}(\M)$ (\cite{DGT}) and so, the
equality
$$c_{\mathcal{R}(\A)}^{\mbox{\scriptsize reg}}(m)\cdot (d\cdot d^\prime) = 
c_{\mathcal{R}(\A)}^{\mbox{\scriptsize ort}}(m)\cdot d^\prime \cdot 1_X$$
implies the existence of a morphism $t$ such that
$t\cdot d\cdot d^\prime=d^\prime$. 
Consequently, $d$ is an isomorphism and 
$c_{\mathcal{R}(\A)}^{\mbox{\scriptsize reg}}(m) \cong
c_{\mathcal{R}(\A)}^{\mbox{\scriptsize ort}}(m)$.
\end{proof}

\section{Injectivity and categories of algebras}

Let $A$, $\mathcal{R}(A)$ and $\mathcal{E}(A)$ be a topological space, the
reflective hull of $A$ and the epireflective hull of $A$ in $\T$,
respectively. 
It is well-known that the functor 
$U=\mbox{Hom}(-,A):\T^{\op} \rightarrow \S$ is a right adjoint whose left
adjoint is given by the power functor $A^-$. 
Let ${K}^\prime$ be the corresponding comparison functor and let $R$ be
the reflection functor from $\T$ to $\mathcal{R}(A)$. 
It is known that the restriction of $U$ to $\mathcal{R}(A)$ is also a
right adjoint and induces the same monad as $U$.
Furthermore, if $K$ is the comparison functor concerning to the last
adjunction, then, up to isomorphism, 
${K}^\prime=K \cdot R^{\mbox{\scriptsize op}}$, 
as illustrated by the following diagram.
$$
\begin{diagram}
\T^{\mbox{op}}	&\rTo_{R^{\mbox{op}}}		&\mathcal{R}(A)^{\mbox{op}} 	&\rTo_{K}	&\S^{\bkT}	\\
		&\rdTo_{{U=\mbox{Hom}(-,A)}}	&\dTo{}				&\ldTo_{U^\bkT}	&		\\
		&				&\S				&		&		\\
\end{diagram}
$$

We are going to show that, under a convenient injectivity assumption on
$A$, $\mathcal{R}(A)^{\op}$ is a reflective subcategory of $\S^\bkT$. 
We also give conditions under which $K$ is an equivalence.

We are going to make use of the following result:

\begin{smt}[G. Janel\-idze, private communication]
Let $\A$ and $\X$ be categories with coequalizers. 
A functor $U : \A \rightarrow \X$ is monadic if the following conditions
hold:
\begin{itemize}
\item[(a)] $U$ has a left adjoint $F$;
\item[(b)] $U$ reflects isomorphisms;
\item[(c)] the counit $FU \rightarrow 1_\A$ is a split epimorphism.
\end{itemize}
\end{smt}

\begin{teo}\label{theo}
Let $A$ be a topological space such that epimorphisms in $\mathcal{E}(A)$
are left-cancellable with respect to $\M$ and 
$c_{\mathcal{R}(A)}^{\mbox{\scriptsize reg}}$ is weakly hereditary in
$\mathcal{E}(A)$. 
Then, if 
$P_{\mathcal{E}(A)}(\M)=A^{\mbox{\scriptsize Inj}_{\mathcal{E}(A)}}$, 
the comparison functor $K$ is full and faithful and 
$\mathcal{R}(A)^{\mbox{\scriptsize op}}$ is equivalent to a reflective
subcategory of the Eilenberg-Moore category $\S^{\bkT}$.

Furthermore, if 
$P_{\mathcal{E}(A)}(\M) = 
(\mathcal{R}(A))^{\mbox{\scriptsize Inj}_{\mathcal{E}(A)}}$, 
then $K$ is an equivalence and so $\mathcal{R}(A)^{\mbox{\scriptsize op}}$
is monadic.
\end{teo}

\begin{proof}
Under the stated conditions, Proposition \ref{r=c} ensures that the
closure operators $c^{\mbox{\scriptsize reg}}_{\mathcal{R}(\A)}$ and
$c^{\mbox{\scriptsize ort}}_{\mathcal{R}(\A)}$ coincide in 
$P_{\mathcal{E}(\A)}$.

Since Hom$(-,A)$ is a right adjoint and 
$\mathcal{R}(\A)^{\mbox{\scriptsize op}}$ has coequalizers (because it is
a full reflective subcategory of $\T^{\mbox{\scriptsize op}}$), we know
that $K$ is a right adjoint. 
In order to show that $K$ is full and faithful, it suffices to prove that
the co-units of the adjunction Hom$(-,A)$ are regular monomorphisms in
$\mathcal{R}(A)$.
For each $B\in \mathcal{R}(\A)$, the co-unit is given by
$\varepsilon_B:B\rightarrow A^{\mbox{\scriptsize Hom}(B,A)}$,
defined as the unique morphism such that 
$\pi_g\cdot \varepsilon_B=g$ for each $g\in \mbox{Hom}(B,A)$, where
$\pi_g$ is the corresponding projection. 
Since $\varepsilon_B$ clearly belongs to 
$\A^{\mbox{\scriptsize Inj}_{\mathcal{E}(\A)}}$, it belongs to
$P_{\mathcal{E}(\A)}(\M)$. 
But, by Theorem \ref{abs clo}, $B$ is an absolutely 
$c^{\mbox{\scriptsize ort}}_{\mathcal{R}(\A)}$-closed object. 
Consequently, $\epsilon_B$ is 
$c^{\mbox{\scriptsize ort}}_{\mathcal{R}(\A)}$-closed and, since
$c^{\mbox{\scriptsize reg}}_{\mathcal{R}(\A)}$ and
$c^{\mbox{\scriptsize ort}}_{\mathcal{R}(\A)}$ coincide in 
$P_{\mathcal{E}(\A)}$, $\varepsilon_B$ is 
$c^{\mbox{\scriptsize reg}}_{\mathcal{R}(\A)}$-closed, so it is a regular
monomorphism.

In order to prove the second part of the theorem, let
$$P_{\mathcal{E}(A)}(\M) = {\mathcal{R}(A)}^{\Inj_{\mathcal{E}(A)}}.$$
Then each $P_{\mathcal{E}(A)}(\M)$-morphism with domain in
${\mathcal{R}(A)}$ is a split monomorphism. 
In particular, for each $B\in {\mathcal{R}(A)}$, the morphism
$\varepsilon_B:B\rightarrow A^{\mbox{\scriptsize Hom}(B,A)}$ is a split
monomorphism. 
Thus, using the Split Monadicity Theorem, to conclude that 
$\mbox{Hom}(-,A):\mathcal{R}(A)^{\mbox{\scriptsize op}}\rightarrow \S$ 
is monadic, it remains to show that it reflects isomorphisms.
Let $f:B\rightarrow C$ be such that
$\mbox{Hom}(f,A):\mbox{Hom}(C,A)\rightarrow \mbox{Hom}(B,A)$ 
is a bijective function. 
This means that $f\in A^\bot$. 
But
$$
A^\bot = \left(\mathcal{R}(A)\right)^\bot 
\subseteq 
{\mathcal{R}(A)}^{\Inj_\mathcal{E}(A)}\cap \mbox{Epi}(\mathcal{R}(A)).
$$
Thus $f$ is a split monomorphism, because it is a morphism of 
${\mathcal{R}(A)}^{\Inj _\mathcal{E}(A)}$ with domain in
$\mathcal{R}(A)$, and it is an epimorphism in $\mathcal{R}(A)$; therefore
it is an isomorphism.
\end{proof}

\begin{obs}
We point out that the equality
$P_{\mathcal{E}(A)}(\M)={\mathcal{R}(A)}^{\Inj _{\mathcal{E}(A)}}$ 
implies that the reflective hull of $A$ is just the subcategory of all
retracts of powers of $A$. 
Furthermore, the converse implication is true, whenever 
$P_{\mathcal{E}(A)}(\M)=A^{\Inj_{\mathcal{E}(A)}}$.
\end{obs}

As an immediate consequence of Theorem \ref{theo}, we get Corollary below,
which was proved by M. Sobral in \cite{So}. 
We need to make use of the following lemma.

\begin{lem}\label{lf}
If $\A$ is a subcategory of $\T$ such that 
$\A^{\mbox{\scriptsize Inj}_{\mathcal{E}(A)}}$ coincides with the class
of embeddings in $\mathcal{E}(\A)$, then, in $\mathcal{E}(\A)$, embeddings
are pushout-stable and $c_{\mathcal{R}(\A)}^{\reg}$ is weakly hereditary.
\end{lem}

\begin{proof}
It is well-known that each object in $\mathcal{E}(\A)$ is the domain of
some initial monosource with codomain in $\A$ (see \cite{AHS}). 
This fact, combined with the assumption of the lemma on $\A$, is easily
seen to imply that embeddings are pushout-stable in $\mathcal{E}(\A)$.

In order to show that $c_{\mathcal{R}(\A)}^{\reg}$ (which coincides with
$c_{\mathcal{E}(\A)}^{\reg}$) is weakly hereditary in $\mathcal{E}(\A)$,
let $m:X\rightarrow Y$ be an embedding in $\mathcal{E} (\A)$ and let
$(u,v)$ be the cokernel pair of $m$ in $\mathcal{E}(\A)$. 
Then 
$\overline{m}=c_{\mathcal{R}(\A)}^{\mbox{\reg}}(m)=\mbox{eq}(u,v)$. 
Let $d$ be the morphism such that $m=\overline{m}\cdot d$. 
We want to show that $d$ is an epimorphism. 
Let $f$ and $g$ be morphisms such that $fd=gd$. 
As $\mathcal{E}(\A)$ is the epireflective hull of $\A$, we may assume, 
without loss of generality, that the codomain of $f$ and $g$ belongs to
$\A$. 
Since $\overline{m}\in \A^{\mbox{\Inj}(\mathcal{E}(\A))}$, there are 
morphisms $\overline{f}$ and $\overline{g}$ such that 
$\overline{f}\cdot \overline{m}=f$ and $\overline{g}\cdot \overline{m}=g$. 
Then, we get 
$\overline{f}\cdot m = \overline{f}\cdot \overline{m}\cdot d
= \overline{g}\cdot \overline{m}\cdot d
= \overline{g}\cdot m$, 
what
implies the existence of a morphism $t$ such that $\overline{f}=t\cdot u$
and $\overline{g}=t\cdot v$. 
Consequently,
$f = \overline{f}\cdot \overline{m}
= t\cdot u\cdot \overline{m}
= t\cdot v\cdot \overline{m}
= \overline{g}\cdot \overline{m}
= g$.
\end{proof}

\begin{cor}[\cite{So}]
Let $A$ be an $\M$-injective topological space. 
Then the comparison functor $K$ is full and faithful and 
$\mathcal{R}(A)^{\mbox{\scriptsize op}}$ is equivalent to a reflective
subcategory of the Eilenberg-Moore category $\S^{\bkT}$.
\end{cor}

\begin{proof}
The fact that $A$ is $\M$-injective implies that
$\M=A^{\mbox{\scriptsize Inj}_{\mathcal{E}(A)}}$. 
Thus, by Lemma \ref{lf}, $P_{\mathcal{E}(A)}(\M)$ is just the class of all 
embeddings in $\mathcal{E}(A)$, and, from Lemma \ref{M}, this implies that
epimorphisms in $\mathcal{E}(A)$ are left-cancellable with respect to 
$\M$. 
Furthermore, Lemma \ref{lf} assures that 
$c_{\mathcal{R}(A)}^{\mbox{\scriptsize reg}}$ is weakly hereditary. 
Therefore, we are under the conditions of Theorem \ref{theo}.
\end{proof}

\begin{exs}
\mbox{}
\begin{enumerate}
\item
The reflective hull of the topological space 
$$A=( \{0,1,2\}, \{ \emptyset, \{0\}, \{1,2\}, \{0,1,2\} \} )$$ 
is $\T$ (cf.\ \cite{Ho}). 
Since in $\T$ embeddings are pushout-stable and $A$ is $\M$-injective, it
follows, from Theorem \ref{theo}, that the dual of $\T$ is a full
reflective subcategory of the Eilenberg-Moore category induced by the
right adjoint Hom$(-,A):\T^{\op}\rightarrow \S$.
\item
Analogously for $\T_0$, we have that embeddings are pushout-stable and the
Sierpi\'nski space $S$ is $\M$-injective, so that, since the category
$\mathcal{S}ob$ of sober spaces is the reflective hull of $S$, the dual
category of $\mathcal{S}ob$ is a reflective subcategory of the
Eilenberg-Moore category induced by the right adjoint Hom$(-,S)$.
\item
In the category 0-$\T_2$ of 0-dimensional Hausdorff spaces, epimorphisms
are just the dense continuous maps, that is, the continuous maps
$f:X\rightarrow Y$ such that, for each $y$ and each clopen $C$ in $Y$, it
holds that $y\in C\Rightarrow C\cap f(X)\not=\emptyset$. 
Thus, it is clear that, in 0-$\T_2$, epimorphisms are left-cancellable
with respect to $\M$. 
On the other hand, given a subspace $X$ of $Y$ in 0-$\T_2$, being 
$m:X\rightarrow Y$ the corresponding embedding, it is easily seen that the
regular closure of $X$ in $Y$ is just the intersection of all clopens
which contain $X$, and that it is just the Kuratowski closure of $X$ in
$Y$. 
Therefore the morphism, which composed with $c^{\reg}(m)$ is equal to $m$,
is dense, so an epimorphism; that is, $c^{\reg}$ is weakly hereditary. As
seen in \ref{Ex1.1}.4, $P_{}=D^{\Inj}$. 
The reflective hull of $D$ is the category of 0-dimensional compact
Hausdorff spaces (cf.\ \cite{Ho}). 
Therefore, applying the corollary of Theorem \ref{theo}, we conclude that
the opposite category of 0-dimensional compact Hausdorff spaces is a full
reflective subcategory of an Eilenberg-Moore category.
\item
Another topological space which fulfils the conditions of Theorem 
\ref{theo} is the unit interval $I$. 
The reflective hull of $I$ is the category of compact Hausdorff spaces,
and its epireflective hull is the category $\mathcal{T}ych$ of Tychonoff
spaces. 
As mentioned in \ref{Ex1.1}.5, 
$P_{\mathcal{T}ych}(\M)=I^{\Inj _{\mathcal{T}ych}}$. 
The regular closure operator $c^{\reg}_{\mathcal{R}(I)}$ in 
${\mathcal{T}ych}$ coincides with the Kuratowski closure
(cf.\ \cite{DT}); consequently it is weakly hereditary, and epimorphisms 
are dense morphisms, which are left-cancellabble with respect to $\M$.
\item
The category $\mathcal{I}nd$ of indiscrete spaces is the epireflective
hull and the reflective hull of the indiscrete space $\{0,1\}$
(cf.\ \cite{Ho}). 
It holds 
$P_{\mathcal{I}nd}(\M) = \M 
= \{ 0,1\}^{\Inj_\mathcal{I}nd} 
= {\mathcal{R}(\{0,\\})^{\Inj_\mathcal{I}nd}$.
Then $\mathcal{I}nd^{\op}$ is equivalent to the Eilenberg-Moore category
induced by the two-points indiscrete space.
\end{enumerate}
\end{exs}

\providecommand{\bysame}{\leavevmode\hbox to3em{\hrulefill}\thinspace}
\providecommand{\MR}{\relax\ifhmode\unskip\space\fi MR }
\providecommand{\MRhref}[2]{%
  \href{http://www.ams.org/mathscinet-getitem?mr=#1}{#2}
}
\providecommand{\href}[2]{#2}

\end{document}